# CRAMÉR-TYPE LARGE DEVIATIONS FOR SAMPLES FROM A FINITE POPULATION[1]

BY ZHISHUI HU, JOHN ROBINSON AND QIYING WANG

*USTC, University of Sydney and University of Sydney*

Cramér-type large deviations for means of samples from a finite population are established under weak conditions. The results are comparable to results for the so-called self-normalized large deviation for independent random variables. Cramér-type large deviations for the finite population Student $t$-statistic are also investigated.

**1. Introduction and results.** Let $X_1, X_2, \ldots, X_n$ be a simple random sample drawn without replacement from a finite population $\{a\}_N = \{a_1, \ldots, a_N\}$, where $n < N$. Denote $\mu = EX_1$, $\sigma^2 = \text{var}(X_1)$,

$$S_n = \sum_{k=1}^n X_k, \qquad p = n/N, \qquad q = 1 - p, \qquad \omega_N^2 = Npq.$$

Under appropriate conditions, the finite central limit theorem (see [14]) states that $P(S_n - n\mu \geq x\sigma\omega_N)$ may be approximated by $1 - \Phi(x)$, where $\Phi(x)$ is the distribution function of a standard normal variate. The absolute error of this normal approximation, via Berry–Esseen bounds and Edgeworth expansions, has been widely investigated in the literature. We only refer to [4] and [17] for the rates in the Erdös and Rényi central limit theorem and to [1, 3, 26], as well as [6, 7], for the Edgeworth expansions. Extensions to $U$-statistics and, more generally, symmetric statistics can be found in [20, 21, 27, 28], as well as [9, 10].

In this paper, we shall be concerned with the relative error of $P(S_n - n\mu \geq x\sigma\omega_N)$ to $1 - \Phi(x)$. In this direction, Robinson [25] derived a large deviation result that is similar to the type for sums of independent random variables in [22], Chapter VIII. However, to make the main results in [25] applicable, it essentially requires the assumption that $0 < p_1 \leq p \leq p_2 < 1$. This kind of

Received November 2004; revised May 2006.
[1]Supported in part by ARC DP0451722.

condition not only takes away a major difficulty in proving large deviation results but also limits its potential applications. The aim of this paper is to establish a Cramér-type large deviation for samples from a finite population under weak conditions. In a reasonably wide range for $x$, we show that the relative error of $P(S_n - n\mu \geq x\sigma\omega_N)$ to $1 - \Phi(x)$ is only related to $E|X_1 - \mu|^3/\sigma^3$ by means of an absolute constant. We also obtain a similar result for the so-called finite population Student $t$-statistic defined by

$$t_n = \sqrt{n}(\bar{X} - \mu)/(\hat{\sigma}\sqrt{q}),$$

where $\bar{X} = S_n/n$ and $\hat{\sigma}^2 = \sum_{j=1}^n (X_j - \bar{X})^2/(n-1)$. It is interesting to note that the results for both the finite population standardized mean and the Student $t$-statistic are comparable to the so-called self-normalized large deviation for independent random variables which has recently been developed by Jing, Shao and Wang [19]. Indeed, Theorems 1.1 and 1.3 below can be considered as analogous to Theorem 2.1 of Jing, Shao and Wang [19] in the independent case. The Berry–Esseen bounds and Edgeworth expansions for the Student $t$-statistic have been investigated in [2], [23] and [5, 8].

We now state our main findings.

THEOREM 1.1. *There exists an absolute constant $A > 0$ such that*

$$\exp\{-A(1+x)^3\beta_{3N}/\omega_N\} \leq \frac{P(S_n - n\mu \geq x\sigma\omega_N)}{1 - \Phi(x)}$$
(1.1)
$$\leq \exp\{A(1+x)^3\beta_{3N}/\omega_N\}$$

*for $0 \leq x \leq (1/A)\omega_N\sigma/\max_k |a_k - \mu|$, where $\beta_{3N} = \sigma^{-3}E|X_1 - \mu|^3$.*

The following result is a direct consequence of Theorem 1.1 and provides a Cramér-type large deviation result for samples from a finite population.

THEOREM 1.2. *There exists an absolute constant $A > 0$ such that*

(1.2) $$\frac{P(S_n - n\mu \geq x\sigma\omega_N)}{1 - \Phi(x)} = 1 + O(1)(1+x)^3\beta_{3N}/\omega_N$$

*and*

(1.3) $$\frac{P(S_n - n\mu \leq -x\sigma\omega_N)}{\Phi(-x)} = 1 + O(1)(1+x)^3\beta_{3N}/\omega_N$$

*for $0 \leq x \leq (1/A)\min\{\omega_N\sigma/\max_k |a_k - \mu|, (\omega_N/\beta_{3N})^{1/3}\}$, where $O(1)$ is bounded by an absolute constant. In particular, if $\omega_N/\beta_{3N} \to \infty$, then, for any $0 < \eta_N \to 0$,*

(1.4) $$\frac{P(S_n - n\mu \geq x\sigma\omega_N)}{1 - \Phi(x)} \to 1, \qquad \frac{P(S_n - n\mu \leq -x\sigma\omega_N)}{\Phi(-x)} \to 1$$

*uniformly in $0 \leq x \leq \eta_N \min\{\omega_N\sigma/\max_k |a_k - \mu|, (\omega_N/\beta_{3N})^{1/3}\}$.*

Results (1.2) and (1.3) are useful because they provide not only the relative error, but also a Berry–Esseen rate of convergence. Indeed, by the fact that $1 - \Phi(x) \leq 2e^{-x^2/2}/(1+x)$ for $x \geq 0$, we may obtain

$$|P(S_n - n\mu \leq x\sigma\omega_N) - \Phi(x)| \leq A(1+|x|)^2 e^{-x^2/2}\beta_{3N}/\omega_N$$

for $|x| \leq (1/A)\min\{\omega_N\sigma/\max_k |a_k - \mu|, (\omega_N/\beta_{3N})^{1/3}\}$. This provides an exponential nonuniform Berry–Esseen bound for samples from a finite population.

REMARK 1.1. We do not have any restriction on the $\{a\}_N$ in Theorems 1.1 and 1.2. Indeed, for any $\{a\}_N$,

$$\mu = \frac{1}{N}\sum_{k=1}^{N} a_k, \qquad \sigma^2 = \frac{1}{N}\sum_{k=1}^{N}(a_k - \mu)^2, \qquad E|X_1 - \mu|^3 = \frac{1}{N}\sum_{k=1}^{N}|a_k - \mu|^3.$$

Removing the trivial case that all $a_k$ are the same, we always have $\max_k |a_k - \mu| > 0$, $\sigma^2 > 0$ and $E|X_1 - \mu|^3 < \infty$.

REMARK 1.2. Hájek [15] proved that if $0 < p_1 \leq p \leq p_2 < 1$, then $(S_n - n\mu)/\sigma\omega_N \to_{\mathcal{D}} N(0,1)$ if and only if $\omega_N\sigma/\max_k |a_k - \mu| \to \infty$. Theorems 1.1 and 1.2 therefore provide reasonably wide ranges for $x$ to make the results hold. To be more precise, as an example consider $a_k = k^\alpha$, where $\alpha > -1/3$. In this special case, simple calculations show that

$$\min\left\{\omega_N\sigma/\max_k |a_k - \mu|, (\omega_N/\beta_{3N})^{1/3}\right\} \asymp (Npq)^{1/6}.$$

Thus, Theorem 1.2 holds for $x$ belonging to the best range $(0, o[(Npq)^{1/6}])$.

The following theorem provides a relative error of $P(t_n \geq x)$ to $1 - \Phi(x)$, which is only related to $E|X_1 - \mu|^3/\sigma^3$ by means of an absolute constant, as in Theorem 1.1. Cramér-type large deviation results for the Student $t$-statistic may be obtained accordingly, as in Theorem 1.2. We omit the details.

THEOREM 1.3. *There exists an absolute constant $A > 0$ such that*

(1.5)
$$\exp\{-A(1+x)^3\beta_{3N}/\omega_N\} \leq \frac{P(t_n \geq x)}{1 - \Phi(x)}$$
$$\leq \exp\{A(1+x)^3\beta_{3N}/\omega_N\}$$

*for all $0 \leq x \leq (1/A)\omega_N\sigma/\max_k |a_k - \mu|$, where $\beta_{3N}$ is defined as in Theorem 1.1.*

TABLE 1
*Relative errors of tail probability to normal*

| Population | $N$ | $Np$ | $\beta_{3N}$ | $x=2$ | $x=2.5$ | $x=3$ |
|---|---|---|---|---|---|---|
| $a_k = k$ | 1000 | 250 | 1.299 | 1.006 | 1.041 | 1.119 |
|  | 100 | 25 | 1.299 | 1.223 | 1.562 | 2.356 |
| $a_k = k^2$ | 1000 | 250 | 1.376 | 0.897 | 0.807 | 0.707 |
|  | 100 | 25 | 1.375 | 0.802 | 0.783 | 0.814 |

REMARK 1.3. The finite population $\{a\}_N$ in Theorems 1.1–1.3 may be replaced by a triangular array $\{a\}_N = \{a_{Ni}, i = 1, 2, \ldots, N\}, N \geq 2$, without essential difficulty. Theorems 1.1–1.2 give information on normal approximations for permutation and rank tests, while Theorem 1.3 is applicable in survey sampling. It is worthwhile to note that we did not introduce any restrictions on unknown parameters or the sampled population and, by means of methods similar to those used in the proof of Theorem 1.3, it is possible to obtain similar bounds for the Studentized mean under stratified random sampling.

REMARK 1.4. The importance of Theorems 1.1–1.3 is based on the fact that all bounds are based only on $E|X_1 - \mu|^3/\sigma^3$ with an absolute constant. The relevance of the results lies in the fact that they give general bounds on the relative errors of the normal approximation, in the same way as the Berry–Esseen bounds are of use in giving uniform bounds on the absolute error. However, our theorems are still asymptotic results, as the absolute constant $A$ is not specified. The following simulations in Table 1, which provide the relative error of $P(t_n \geq x)$ to $1 - \Phi(x)$ based on Theorem 1.3 for $a_k = k$ and $a_k = k^2$ based on 1,000,000 repetitions, indicate the accuracy of the normal approximation in the large deviation region. The simulations in Table 2 confirm that more precise results, like saddle-point approximations, are required. Dai and Robinson [11] investigated saddle-point approximation under very strong conditions, that is, where errors are based on $Ee^{t(X_1-\mu)^2}$, for $t > 0$, with a constant. The saddle-point approximations under the weak conditions as in Theorem 1.3 are open problems.

This paper is organized as follows. Major steps of the proofs of Theorems 1.1–1.3 are given in Section 2. Proofs of three propositions used in the main proofs are offered in Sections 3–5. In Section 3, as a preliminary, we establish a Berry–Esseen bound for the associated distribution of $P(S_n - n\mu \leq x)$ related to the conjugate method in a general setting. Throughout the paper, we shall use $A, A_1, A_2, \ldots$ to denote absolute constants whose values may differ at each occurrence. We also write $b = x/\omega_N$, $V_n^2 = \sum_{k=1}^n X_k^2$,

$$V_{1n} = V_n^2 - n \quad \text{and} \quad V_{2n} = \sum_{k=1}^n [(X_k^2 - 1)^2 - E(X_k^2 - 1)^2]$$

and, when no confusion arises, $\sum$ denotes $\sum_{k=1}^N$ and $\prod$ denotes $\prod_{k=1}^N$. The symbol $i$ will be used exclusively for $\sqrt{-1}$.

**2. Proofs of theorems.** Without loss of generality, we assume that $\mu = 0$ and $\sigma^2 = 1$. Otherwise, it suffices to consider that $\{X_1, X_2, \ldots, X_n\}$ is a simple random sample drawn without replacement from a finite population $\{a'\}_N = \{(a_1 - \mu)/\sigma, \ldots, (a_N - \mu)/\sigma\}$, where $n < N$.

PROOF OF THEOREM 1.1. When $0 \leq x \leq 2$, property (1.1) follows from the Berry–Esseen bound for samples from a finite population (see, e.g., [17]),

$$|P(S_n \geq x\omega_N) - (1 - \Phi(x))| \leq A\beta_{3N}/\omega_N.$$

When $2 \leq x \leq (1/A)\omega_N/\max_k |a_k|$, property (1.1) follows from the following proposition with $\xi = 0$, $\xi_1 = 0$ and $h = 0$. Proposition 2.1 will be proven in Section 3. □

PROPOSITION 2.1. *There exists an absolute constant $A > 0$ such that, for all $0 \leq \xi \leq 1/2$, $|\xi_1| \leq 36$ and $2 \leq x \leq (1/A)\omega_N/\max_k |a_k|$,*

$$(2.1) \quad \frac{P(bS_n - \xi b^2 qV_{1n} + \xi_1 b^4 q^2 V_{2n} \geq x^2)}{1 - \Phi(x)} \geq \exp\{-Ax^3\beta_{3N}/\omega_N\}$$

*and*

$$(2.2) \quad \frac{P(bS_n - \xi b^2 qV_{1n} + \xi_1 b^4 q^2 V_{2n} \geq x^2 + h)}{1 - \Phi(x)} \leq [1 + 9|h|x^{-2}]\exp\{-h + Ax^3\beta_{3N}/\omega_N\},$$

*where $h$ is an arbitrary constant (which may depend on $x$) with $|h| \leq x^2/5$.*

TABLE 2
*Saddle-point approximations*

| Population | $N$ | $Np$ | $\beta_{3N}$ | $x = 2$ | $x = 2.5$ | $x = 3$ |
|---|---|---|---|---|---|---|
| $a_k = k$ | 1000 | 250 | 1.299 | 0.984 | 0.989 | 1.011 |
| | 100 | 25 | 1.299 | 0.996 | 0.989 | 0.992 |
| $a_k = k^2$ | 1000 | 250 | 1.376 | 1.010 | 0.987 | 0.966 |
| | 100 | 25 | 1.375 | 0.991 | 1.007 | 1.009 |

PROOF OF THEOREM 1.2.  This follows immediately from Theorem 1.1. □

PROOF OF THEOREM 1.3.  When $0 \leq x \leq 4$, property (1.5) follows from the Berry–Esseen bound for the finite population Student $t$-statistic; see, for example, [5]. Next, assume $4 \leq x \leq (1/A)\omega_N/\max_k |a_k|$. Without loss of generality, assume that $A \geq 8$ and $n \geq 4$. Note that $\max_k |a_k| \geq 1$, since $\sum a_k^2 = N$. It is readily seen that

$$(2.3) \qquad \left|\frac{x_0}{x} - 1\right| = |[1 + (x^2 q - 1)/n]^{-1/2} - 1| \leq 2x^2/n,$$

where $x_0 = xn^{1/2}/(n + x^2 q - 1)^{1/2}$. It follows from (2.3) that $2 \leq x/2 \leq x_0 \leq 3x/2$ and $|x_0 - x| \leq 2x^3 \beta_{3N}/\omega_N^2$. Hence, by noting $1 - \Phi(x) \geq x\Phi'(x)/(1+x^2)$ for $x \geq 0$ (see, e.g., [24], page 30), we have

$$\left|\log \frac{1 - \Phi(x_0)}{1 - \Phi(x)}\right| = \left|\int_x^{x_0} \frac{\Phi'(t)}{1 - \Phi(t)} dt\right| \leq \left|\int_x^{x_0} \frac{1 + t^2}{t} dt\right|$$
$$\leq 2x|x - x_0| \leq x^3 \beta_{3N}/\omega_N,$$

which yields that

$$(2.4) \qquad \exp\{-x^3 \beta_{3N}/\omega_N\} \leq \frac{1 - \Phi(x_0)}{1 - \Phi(x)} \leq \exp\{x^3 \beta_{3N}/\omega_N\}.$$

We are now ready to prove Theorem 1.3. As is well known (see, e.g., [13]), for $x \geq 0$,

$$P(t_n \geq x) = P(S_n/V_n \geq x_0 \sqrt{q}).$$

Note that $b_0 x_0 \sqrt{q} V_n \leq (x_0^2 + b_0^2 q V_n^2)/2 \leq x_0^2 + b_0^2 q(V_n^2 - n)/2$, where $b_0 = x_0/\omega_N$. It follows from (2.1), (2.3) and (2.4) that, for $4 \leq x \leq (1/A)\omega_N/\max_k |a_k|$,

$$P(S_n \geq x_0\sqrt{q}V_n) \geq P(b_0 S_n - b_0^2 q(V_n^2 - n)/2 \geq x_0^2)$$
$$\geq (1 - \Phi(x_0))\exp\{-Ax_0^3 \beta_{3N}/\omega_N\}$$
$$\geq (1 - \Phi(x))\exp\{-A_1 x^3 \beta_{3N}/\omega_N\},$$

which implies the first inequality of (1.5).

In view of the following two propositions, the second inequality of (1.5) may be obtained by an argument similar to that used in the proof of (5.13) in [19], so the details are omitted. The proofs of Propositions 2.2 and 2.3 will be given in Section 4 and Section 5, respectively. □

PROPOSITION 2.2. *There exists an absolute constant $A > 0$ such that*

$$P(S_n \geq x\sqrt{q}V_n) \leq (1 - \Phi(x))\exp\{Ax^3 \beta_{3N}/\omega_N\} + Ae^{-4x^2}$$

*for $2 \leq x \leq (1/A)\omega_N/\max_k |a_k|$.*

PROPOSITION 2.3. *There exists an absolute constant $A > 0$ such that*

$$P(S_n \geq x\sqrt{q}V_n) \leq (1 - \Phi(x))\exp\{Ax^3\beta_{3N}/\omega_N\} + A(x\beta_{3N}/\omega_N)^{4/3}$$

*for $2 \leq x \leq (1/A)\omega_N/\max_k |a_k|$.*

**3. Proof of Proposition 2.1.** In Section 3.1, we derive a Berry–Esseen bound for the associated distribution of $P(S_n \leq x)$ related to the conjugate method. The result is established in a general setting and will be used in the proof of Proposition 2.1, which is given in Section 3.2.

3.1. *Preliminary.* Consider a sequence of constants $\{b\}_N = \{b_1, \ldots, b_N\}$ with $\sum b_k = 0$. Write $K(z) = \log(pe^{qz} + qe^{-pz})$ and let $K_k, K'_k$ and $K''_k$ be the values of $K(x), K'(x)$ and $K''(x)$ evaluated at $x = ub_k + \alpha_N$, where $\alpha_N$ is the solution of the equation

(3.1) $$\sum K'(ub_k + \alpha) = 0.$$

We continue to assume that $X_1, X_2, \ldots, X_n$ is a random sample without replacement from $\{b\}_N$, where $n < N$, and continue to use the notation $S_n = \sum_{k=1}^n X_k$, $p, q$ and $\omega_N^2 = Npq$, as in Section 1. Define

$$H_n(x; u) = Ee^{uS_n}I(S_n \leq x)/Ee^{uS_n}.$$

The main result in this section is as follows.

THEOREM 3.1. *For any given $C_0 > 0$, if $|u| \leq C_0/\max_k |b_k|$, then*

(3.2) $$\sup_x \left| H_n(x; u) - \Phi\left(\frac{x - m_N}{\sigma_N}\right) \right| \leq C(pq)^{-1/2}\sum |b_k|^3 / \left(\sum b_k^2\right)^{3/2}$$

*and*

(3.3) $$Ee^{uS_n} = (G_n(p))^{-1}\left(\sum K''_k\right)^{-1/2}\exp\left\{\sum K_k\right\}(1 + O_1/\omega_N),$$

*where $G_n(p) = \sqrt{2\pi}\binom{N}{n}p^n q^{N-n}$, $m_N = \sum b_k K'_k$,*

$$\sigma_N^2 = \sum b_k^2 K''_k - \left(\sum b_k K''_k\right)^2 / \sum K''_k,$$

*$|O_1|$ is bounded by $C_1$ and both $C$ and $C_1$ are constants depending only on $C_0$.*

The results (3.2) and (3.3) essentially improve Lemma 2 and Lemma 1 (with $v = 0$) of [25], respectively. In [25], the constants $C$ and $C_1$ are imposed to depend on $p$ and $q$. The proof of Theorem 3.1 follows the approach of Robinson [25], but with quite different calculations. The details can be found in [18], on which the present paper is based.

3.2. *Proof.* Roughly speaking, the proof of Proposition 2.1 is based on the conjugate method and an application of Theorem 3.1 to the $b_k$ specified in (3.4) below.

Let $0 < \lambda \leq 2$, $0 \leq \theta \leq 1$ and $|\theta_1| \leq 72$. Define, for $k = 1, \ldots, N$,

$$(3.4) \quad b_k = \lambda b a_k - \theta b^2 q(a_k^2 - 1) + \theta_1 b^4 q^2 \left[ (a_k^2 - 1)^2 - \frac{1}{N} \sum (a_j^2 - 1)^2 \right].$$

Since $\sum a_k = 0$ and $\sum a_k^2 = N$, it is readily seen that $\max_k |a_k| \geq 1$ and $\sum b_k = 0$. Also, when $b \max_k |a_k| \leq 1/128$, we have that $b\beta_{3N} \leq 1/128$,

$$(3.5) \qquad \max_k |b_k| \leq 1/32,$$

$$(3.6) \qquad \left| \sum b_k^2 - \lambda^2 b^2 N \right| \leq 5Nb^3 q \beta_{3N},$$

$$(3.7) \qquad \sum |b_k|^3 \leq 9Nb^3 \beta_{3N}.$$

Recalling $b = x/\omega_N$, (3.5)–(3.7) hold if $0 \leq x \leq (1/128)\omega_N / \max_k |a_k|$.

We establish five lemmas before proceeding to the proof of Proposition 2.1. The first lemma summarizes some basic properties of $K(z)$. We still use the notation $K_k$, $K'_k$ and $K''_k$ to denote the values of $K(z)$, $K'(z)$ and $K''(z)$ evaluated at $z = b_k + \alpha_N$, where the $b_k$ are limited to (3.4) and $\alpha_N$ is the solution of the equation $\sum K'(b_k + \alpha) = 0$. Lemmas 3.2 and 3.3 give the properties of $\alpha_N$, $K_k$, $K'_k$ and $K''_k$. Lemmas 3.4 and 3.5 provide the results that will be used in the proof of Proposition 2.1.

LEMMA 3.1. *We have $K'(0) = 0$,*

$$(3.8) \qquad -pqe^{2t} \leq K'(-x) < 0 < K'(x) \leq pqe^{2t}, \qquad \text{for } 0 < x \leq t,$$

$$(3.9) \qquad pqe^{-3t} < K''(x) < pqe^{3t}, \qquad \text{for } |x| \leq t,$$

$$(3.10) \quad |K'''(x + iy)| \leq 2^{3/2} e^{5t} pq, \qquad \text{for } |x| \leq t \text{ and } |y| \leq \pi/2.$$

*Furthermore, if $|x| \leq 1/16$, then*

$$(3.11) \qquad |K(x)/pq - x^2/2| \leq (1/2)|x|^3,$$

$$(3.12) \qquad |K'(x)/pq - x| \leq x^2,$$

$$(3.13) \qquad |K''(x)/pq - 1 - (q - p)x| \leq 8x^2.$$

PROOF. The proof of Lemma 3.1 is straightforward and the details are omitted. □

LEMMA 3.2. *If $0 \leq x \leq (1/128)\omega_N / \max_k |a_k|$, then $\alpha_N$ is unique,*

$$(3.14) \quad |\alpha_N| \leq \min\left\{1/32, (2/N) \sum b_k^2\right\} \quad \text{and} \quad \alpha_N^2 \leq (9/8) b^3 \beta_{3N}.$$

PROOF. By virtue of (3.5), (3.8) and (3.9), $\sum K'(b_k + \alpha)$ is negative when $\alpha < -1/32$ and positive when $\alpha > 1/32$, and it is strictly monotone in the range $|\alpha| \leq 1/32$. This implies that $\sum K'(b_k + \alpha) = 0$ has a unique solution $\alpha_N$ and that $|\alpha_N| \leq 1/32$. By noting $|b_k| + |\alpha_N| \leq 1/16$, it follows from (3.12), $\sum K'(b_k + \alpha_N) = 0$ and $\sum b_k = 0$ that

$$N|\alpha_N| = \left|\sum [K'(b_k + \alpha_N)/pq - (b_k + \alpha_N)]\right|$$
$$\leq \sum (b_k + \alpha_N)^2 = \sum b_k^2 + N\alpha_N^2$$
$$\leq \sum b_k^2 + N|\alpha_N|/2.$$

This yields $|\alpha_N| \leq (2/N)\sum b_k^2$ and hence the first result of (3.14) follows. Furthermore, by using Hölder's inequality, $|b_k| \leq 1/32$ and (3.7), we obtain

$$\alpha_N^2 \leq (4/N)\sum b_k^4 \leq (9/8)b^3 \beta_{3N},$$

which implies the second result of (3.14). The proof of Lemma 3.2 is complete. □

LEMMA 3.3. *If* $0 \leq x \leq (1/128)\omega_N/\max_k |a_k|$, *then*

(3.15) $$\left|\sum K_k - \lambda^2 x^2/2\right| \leq 24x^3 \beta_{3N}/\omega_N,$$

(3.16) $$\left|\sum b_k K'_k - \lambda^2 x^2\right| \leq 24x^3 \beta_{3N}/\omega_N,$$

(3.17) $$\left|\sum K''_k - \omega_N^2\right| \leq 41x^2,$$

(3.18) $$\left|\sum b_k K''_k\right| \leq 6x^2,$$

(3.19) $$\left|\sum b_k^2 K''_k - \lambda^2 x^2\right| \leq 21x^3 \beta_{3N}/\omega_N.$$

PROOF. We prove (3.15). The other proofs are similar and hence are omitted. Applying (3.11) with $x = b_k + \alpha_N$ and using Hölder's inequality, we obtain

(3.20) $$\left|\sum [K_k - 2^{-1}pq(b_k + \alpha_N)^2]\right| \leq 2pq\left(\sum |b_k|^3 + N|\alpha_N|^3\right).$$

This, together with $\sum b_k = 0$, (3.6)–(3.7) and (3.14), implies that

$$\left|\sum K_k - \lambda^2 x^2/2\right| \leq \left|\sum [K_k - 2^{-1}pq(b_k + \alpha_N)^2]\right|$$
$$+ 2^{-1}pq\left|\sum b_k^2 - \lambda^2 b^2 N\right| + 2^{-1}\omega_N^2 \alpha_N^2$$
$$\leq 24b^3 \omega_N^2 \beta_{3N} = 24x^3 \beta_{3N}/\omega_N,$$

as required. □

Let $Y_j, 1 \leq j \leq n$, be a random sample of size $n$ without replacement from $\{b_1, b_2, \ldots, b_N\}$ defined by (3.4),

$$T_n^* \equiv T_n(\lambda, \theta, \theta_1) = \sum_{k=1}^n Y_k, \qquad m_N^* \equiv m_N(\lambda, \theta, \theta_1) = \sum b_k K_k',$$

$$\sigma_N^{*2} \equiv \sigma_N^2(\lambda, \theta, \theta_1) = \sum b_k^2 K_k'' - \left(\sum b_k K_k''\right)^2 / \sum K_k''$$

and $H_n^*(x) = E\exp(T_n^*)I(T_n^* \leq x)/E\exp(T_n^*)$.

LEMMA 3.4. *There exists an absolute constant $\lambda_0 > 0$ such that, for $2 \leq x \leq \lambda_0 \omega_N / \max_k |a_k|$,*

$$\exp\{\lambda^2 x^2/2 - Ax^3 \beta_{3N}/\omega_N\} \leq E\exp(T_n^*)$$
(3.21)
$$\leq \exp\{\lambda^2 x^2/2 + Ax^3 \beta_{3N}/\omega_N\}.$$

PROOF. Without loss of generality, assume that $\lambda_0 \leq \min\{1/128, 1/(8C_1 + 4)\}$, where $C_1$ is defined as in Theorem 3.1. Recall that $\max_k |b_k| \leq 1/32$, by (3.5). It follows from (3.3) (Theorem 3.1 with $u = 1$ and $C_0 = 1/32$) that

$$(3.22) \quad E\exp(T_n^*) = (G_n(p))^{-1} \left(\sum K_k''\right)^{-1/2} \exp\left\{\sum_{j=1}^N K_k\right\}(1 + R^*),$$

where $G_n(p) = \sqrt{2\pi}\binom{N}{n} p^n q^{N-n}$ and $|R^*| \leq C_1/\omega_N$. By Stirling's formula,

$$\binom{N}{n} p^n q^{N-n} = (2\pi\omega_N^2)^{-1/2}(1 + O_2 \omega_N^{-2}),$$

where $|O_2| \leq 1/6$. This, together with the fact that $\omega_N \geq x \max_k |a_k|/\lambda_0 \geq 128$ (recall that $\max_k |a_k| \geq 1$), implies that

$$(3.23) \qquad \omega_N^{-1} G_n(p)^{-1}(1 + R^*) = 1 + O_3 \omega_N^{-1},$$

where $|O_3| \leq 2C_1 + 1$. On the other hand, it follows from (3.17) that

$$(3.24) \qquad \left(\sum K_k''\right)^{-1/2} \omega_N = 1 + O_4 b^2,$$

where $|O_4| \leq 82$. From (3.23)–(3.24), for $2 \leq x \leq \lambda_0 \omega_N / \max_k |a_k|$,

$$\exp\{-2A_1 x^3 \beta_{3N}/\omega_N\} \leq \left(\sum K_k''\right)^{-1/2} G_n(p)^{-1}(1 + R^*)$$
(3.25)
$$\leq \exp\{A_1 x^3 \beta_{3N}/\omega_N\},$$

where $A_1 = 2C_1 + 83$ and where we have used the fact that $1/\omega_N + b^2 \leq x^3 \beta_{3N}/\omega_N$, since $b = x/\omega_N$ and $\beta_{3N} \geq 1$. (3.21) now follows easily from (3.15), (3.22) and (3.25). The proof of Lemma 3.4 is thus complete. □

LEMMA 3.5. *There exists an absolute constant $\lambda_1 > 0$ such that, for $2 \leq x \leq \lambda_1 \omega_N / \max_k |a_k|$,*

(3.26) $$|m_N^* - \lambda^2 x^2| \leq 24 x^3 \beta_{3N}/\omega_N,$$

(3.27) $$|\sigma_N^{*\,2} - \lambda^2 x^2| \leq 22 x^3 \beta_{3N}/\omega_N.$$

*If, in addition, $1 \leq \lambda \leq 2$, then*

(3.28) $$\Delta_N := \sup_y |H_n^*(u(y)) - \Phi(y)| \leq 12\, C\beta_{3N}/\omega_N \leq 1/4,$$

*where $u(y) = y\sigma_N^* + m_N^*$ and where $C$ is defined as in* (3.2). *Also, for all $y$ satisfying $m_N^* \geq y + 2\sigma_N^*$, we have*

(3.29) $$P(T_n^* \geq y) \geq (1/2) \exp\{-m_N^* - 2\sigma_N^*\} E \exp(T_n^*).$$

PROOF. Without loss of generality, assume that $\lambda_1 \leq \min\{1/128, 1/(25C)\}$, where $C$ is defined as in (3.2). (3.26) and (3.27) then follow from (3.16)–(3.19) by a simple calculation.

If $1 \leq \lambda \leq 2$, then by noting that $\beta_{3N}/\omega_N \leq x\beta_{3N}/(2\omega_N) \leq \min\{1/128, 1/(50C)\}$, since $\beta_{3N} \leq \max_k |a_k|$, it follows easily from (3.5)–(3.7) that $pq \sum b_k^2 \geq 4x^2/5$ and

(3.30) $$(pq)^{-1/2} \sum |b_k|^3 / \left(\sum b_k^2\right)^{3/2} \leq 12\beta_{3N}/\omega_N \leq 1/(4C).$$

By (3.30) and Theorem 3.1 with $C_0 = 1/32$ and $u = 1$ (recall that $\max_k |b_k| \leq 1/32$),

$$\Delta_N \leq C(pq)^{-1/2} \sum |b_k|^3 / \left(\sum b_k^2\right)^{3/2} \leq 12\, C\beta_{3N}/\omega_N \leq 1/4,$$

which implies (3.28).

By (3.28) and the conjugate method, for all $y$ satisfying $m_N^* \geq y + 2\sigma_N^*$,

$$P(T_n^* \geq y)/E \exp(T_n^*) = \int_y^\infty e^{-u}\, dH_n^*(u)$$

$$= e^{-m_N^*} \int_{(y-m_N^*)/\sigma_N^*}^\infty e^{-x\sigma_N^*}\, dH_n^*(u(y))$$

$$\geq e^{-m_N^* - 2\sigma_N^*} \int_{-2}^2 dH_n^*(u(y))$$

$$\geq e^{-m_N^* - 2\sigma_N^*} (P(|N(0,1)| \leq 2) - \Delta_N)$$

$$\geq (1/2) \exp\{-m_N^* - 2\sigma_N^*\},$$

where $N(0,1)$ is a standard normal random variable and where we have used the fact that

$$P(|N(0,1)| \leq 2) > 3/4.$$

This proves (3.29) and also completes the proof of Lemma 3.5. □

After these preliminaries, we are now ready to prove Proposition 2.1.
In addition to the previous notation, we further let $T_{1n} = T_n(1, \xi, \xi_1)$,
$m_{1N} = m_N(1, \xi, \xi_1)$, $\sigma_{1N}^2 = \sigma_N^2(1, \xi, \xi_1)$, $\varepsilon_N = (x^2 + h - m_{1N})/\sigma_{1N}$
and $H_{1n}(t) = E\exp\{T_{1n}\}I(T_{1n} \le t)/E\exp\{T_{1n}\}$. Note that
$$bS_n - \xi b^2 q V_{1n} + \xi_1 b^4 q^2 V_{2n} = T_{1n}.$$

It follows from the conjugate method that

$$
\begin{aligned}
P(bS_n &- \xi b^2 q V_{1n} + \xi_1 b^4 q^2 V_{2n} \ge x^2 + h) \\
&= P(T_{1n} \ge x^2 + h) \\
&= E\exp\{T_{1n}\} \int_{x^2+h}^{\infty} e^{-t} \, dH_{1n}(t) \\
&= E\exp\{T_{1n}\} e^{-x^2-h} \int_0^{\infty} e^{-t\sigma_{1N}} \, dH_{1n}[\sigma_{1N}(t + \varepsilon_N) + m_{1N}] \\
&= E\exp\{T_{1n}\} e^{-x^2-h} (\mathcal{L}_N + R_N),
\end{aligned}
$$
(3.31)

where
$$\mathcal{L}_N = \int_0^{\infty} e^{-t\sigma_{1N}} \, d\Phi(t + \varepsilon_N),$$
$$R_N = \int_0^{\infty} e^{-t\sigma_{1N}} \, d\{H_{1n}[\sigma_{1N}(t + \varepsilon_N) + m_{1N}] - \Phi(t + \varepsilon_N)\}.$$

We next estimate $E\exp\{T_{1n}\}$, $\mathcal{L}_N$ and $R_N$ for $0 \le \xi \le 1/2$, $|\xi_1| \le 36$, $|h| \le x^2/5$ and $2 \le x \le \eta\omega_N/\max_k |a_k|$, where we assume $\eta$ to be sufficiently small so that $\eta \le \min\{1/128, \lambda_0, \lambda_1\}$, with $\lambda_0$ and $\lambda_1$ defined as in Lemmas 3.4 and 3.5. This choice of $\eta$ guarantees that Lemmas 3.2–3.5 hold and, since $\beta_{3N} \le \max_k |a_k|$, we have

(3.32) $\quad\quad\quad\quad \beta_{3N}/\omega_N \le x\beta_{3N}/(2\omega_N) \le \eta/2 \le 1/256.$

Clearly, by Lemma 3.4,

(3.33)
$$\exp\{x^2/2 - Ax^3\beta_{3N}/\omega_N\} \le E\exp\{T_{1n}\}$$
$$\le \exp\{x^2/2 + Ax^3\beta_{3N}/\omega_N\}.$$

In order to estimate $\mathcal{L}_N$, we note that

(3.34)
$$
\begin{aligned}
\mathcal{L}_N &= \frac{1}{\sqrt{2\pi}} \int_0^{\infty} e^{-\sigma_N t - (1/2)(t+\varepsilon_N)^2} \, dt \\
&= \frac{e^{-\varepsilon_N^2/2}}{\sqrt{2\pi}} \int_0^{\infty} e^{-(\varepsilon_N + \sigma_N)t - (1/2)t^2} \, dt \\
&:= \frac{e^{-\varepsilon_N^2/2}}{\sqrt{2\pi}} \mathcal{L}_{1N}.
\end{aligned}
$$

Write $\psi(t) = \{1 - \Phi(t)\}/\Phi'(t) = e^{t^2/2} \int_t^\infty e^{-y^2/2}\, dy$. It is readily seen that

(3.35)
$$3/4 \leq t\psi(t) \leq 1 \qquad (t \geq 2),$$
$$|\psi'(t)| = |t\psi(t) - 1| \leq t^{-2} \qquad (t > 0).$$

On the other hand, $\psi\{\varepsilon_N + \sigma_N\} = \mathcal{L}_{1N}$ and, by virtue of (3.26), (3.27) and (3.32), it follows that

(3.36) $$|\varepsilon_N - h/\sigma_N| \leq 28x^2 \beta_{3N}/\omega_N.$$

If, in addition, we have $|h| \leq x^2/5$, then

(3.37) $$|\varepsilon_N + \sigma_N - x| \leq 3|h|/(2x) + 41x^2 \beta_{3N}/\omega_N \leq 2x/3.$$

Using (3.35)–(3.37), it follows from Taylor's expansion that, for $|h| \leq x^2/5$ and $2 \leq x \leq \eta \omega_N / \max_k |a_k|$,

$$\begin{aligned}
\mathcal{L}_{1N} &= \psi\{\varepsilon_N + \sigma_N\} \\
&= \psi(x) + \psi'(\theta)\{\varepsilon_N + \sigma_N - x\} \qquad [\text{where } \theta \in (x/3, 5x/3)] \\
&= \psi(x)(1 + \tau + O_5 x \beta_{3N}/\omega_N),
\end{aligned}$$

where $|\tau| \leq 9|h|/x^2$ and $|O_5| \leq 120$. Therefore, taking account of (3.34), we obtain for $|h| \leq x^2/5$ and $2 \leq x \leq \eta \omega_N / \max_k |a_k|$,

(3.38) $$\mathcal{L}_N = e^{x^2/2}\{1 - \Phi(x)\}e^{-\varepsilon_N^2/2}(1 + \tau + O_5 x \beta_{3N}/\omega_N).$$

As for $R_N$, by (3.28) and integration by parts,

$$|R_N| \leq 2\sup_t |H_{1n}[\sigma_{1N} t + m_{1N}] - \Phi(t)| \leq 24C\beta_{3N}/\omega_N.$$

This, together with (3.35), implies that for $x \geq 2$,

(3.39) $$R_N = O_6 x \beta_{3N}/\omega_N e^{x^2/2}\{1 - \Phi(x)\},$$

where $|O_6| \leq 32\sqrt{2\pi} C$.

Combining (3.31), (3.33) and (3.38)–(3.39), it is readily seen that for any $|h| \leq x^2/5$ and $2 \leq x \leq \eta \omega_N / \max_k |a_k|$,

$$\frac{P(bS_n - \xi b^2 q V_{1n} + \xi_1 b^4 q^2 V_{2n} \geq x^2 + h)}{1 - \Phi(x)}$$
$$\leq [1 + 9|h|x^{-2}]\exp\{-h + Ax^3 \beta_{3N}/\omega_N\}.$$

This proves (2.2).

Similarly, by letting $h = 0$, it follows from (3.31), (3.33), (3.36), (3.38) and (3.39) that if, in addition, $x^2 \leq \omega_N/\beta_{3N}$, then

$$
\begin{aligned}
\frac{P(bS_n - \xi b^2 q V_{1n} + \xi_1 b^4 q^2 V_{2n} \geq x^2)}{1 - \Phi(x)} & \\
\geq \exp\{-Ax^3 \beta_{3N}/\omega_N - \varepsilon_N^2/2\}[1 - \{|O_5| + |O_6|e^{\varepsilon_N^2/2}\}x\beta_{3N}/\omega_N] & \\
\geq \exp\{-A_1 x^3 \beta_{3N}/\omega_N\}[1 - A_2 x \beta_{3N}/\omega_N] & \\
\geq \exp\{-A_3 x^3 \beta_{3N}/\omega_N\}, &
\end{aligned}
$$
(3.40)

by choosing $\eta$ sufficiently small. From (3.40), property (2.1) will follow if we prove that, for $x^2 \geq \omega_N/\beta_{3N}$ and $2 \leq x \leq \eta \omega_N/\max_k |a_k|$,

(3.41) $$\frac{P(bS_n - \xi b^2 q V_{1n} + \xi_1 b^4 q^2 V_{2n} \geq x^2)}{1 - \Phi(x)} \geq \exp\{-Ax^3 \beta_{3N}/\omega_N\}.$$

We will prove (3.41) by using (3.29). Let $\lambda = 1 + 28x\beta_{3N}/\omega_N$, $\theta = \lambda\xi$ and $\theta_1 = \lambda \xi_1$. Note that $1 \leq \lambda \leq 3/2$ by (3.32), $0 \leq \theta \leq 3/4$ since $0 \leq \xi \leq 1/2$ and $|\theta_1| \leq 72$ since $|\xi_1| \leq 36$. By virtue of (3.26), (3.27), (3.32) and $x^2 \geq \omega_N/\beta_{3N}$, we have $m_N^* \leq \lambda^2 x^2 + 24x^3\beta_{3N}/\omega_N$, $\sigma_N^* \leq 2x \leq 2x^3\beta_{3N}/\omega_N$ and

$$m_N^* - \lambda x^2 - 2\sigma_N^* \geq \lambda(\lambda-1)x^2 - 28x^3\beta_{3N}/\omega_N \geq 0.$$

Now, by (3.29) with $y = \lambda x^2$ and Lemma 3.4, for $x^2 \geq \omega_N/\beta_{3N}$ and $2 \leq x \leq \eta \omega_N/\max_k |a_k|$, it follows that

$$
\begin{aligned}
P(bS_n - \xi b^2 q V_{1n} + \xi_1 b^4 q^2 V_{2n} \geq x^2) &= P(T_n^* \geq \lambda x^2) \\
&\geq \tfrac{1}{2} \exp\{-m_N^* - 2\sigma_N^*\} E \exp\{T_n^*\} \\
&\geq \tfrac{1}{2} \exp\{-x^2/2 - 2x - Ax^3\beta_{3N}/\omega_N\} \\
&\geq (1 - \Phi(x)) \exp\{-A_1 x^3 \beta_{3N}/\omega_N\},
\end{aligned}
$$

which implies (3.41). The proof of Proposition 2.1 is now complete.

**4. Proof of Proposition 2.2.** By the inequality $(1+y)^{1/2} \geq 1 + y/2 - y^2$ for any $y \geq -1$,

$$
\begin{aligned}
P(S_n \geq x\sqrt{q}V_n) & \\
= P\!\left(S_n \geq x\sqrt{nq}\left(1 + \frac{V_n^2 - n}{n}\right)^{1/2}\right) & \\
\leq P\!\left(S_n \geq x\sqrt{nq}\left[1 + \frac{V_{1n}}{2n} - \frac{V_{1n}^2}{n^2}\right]\right) & \\
(4.1) \quad \leq P\!\left(V_{1n}^2 \geq 36x^2\!\left(\sum_{k=1}^n (X_k^2 - 1)^2 + 5p\sum a_k^4\right)\right) &
\end{aligned}
$$

$$+ P\left(S_n \geq x\sqrt{nq}\left(1 + \frac{V_{1n}}{2n} - \frac{36x^2}{n^2}\left(\sum_{k=1}^{n}(X_k^2 - 1)^2 + 5p\sum a_k^4\right)\right)\right)$$

$$:= R_{1n} + R_{2n}, \quad \text{say.}$$

Note that $R_{2n} = P(bS_n - \frac{1}{2}b^2qV_{1n} + 36b^4q^2V_{2n} \geq x^2 - h_0)$, where, whenever $2 \leq x \leq (1/128)\omega_N/\max_k |a_k|$, we have

$$h_0 = \frac{180px^4 \sum a_k^4}{n^2} + \frac{36x^4 \sum_{k=1}^{n} E(X_k^2 - 1)^2}{n^2} \leq \frac{3x^3 \beta_{3N}}{\omega_N}$$

and also $0 \leq h_0 \leq x^2/5$. It follows from Proposition 2.1 with $\xi = 1/2$, $\xi_1 = 36$ and $h = h_0$ that there exists an absolute constant $A > 128$ such that, for all $2 \leq x \leq (1/A)\omega_N/\max_k |a_k|$,

(4.2) $$R_{2n} \leq (1 - \Phi(x))\exp\{Ax^3 \beta_{3N/\omega_N}\}.$$

This, together with (4.1), implies that Proposition 2.2 will follow if we prove for all $x > 0$ that

(4.3) $$R_{1n} \leq 2\sqrt{2}e^{-4x^2}.$$

Theorem 2.1 of [12] will be used to prove (4.3). To use the theorem, let $Y_i = X_i^2 - 1$, $\mathcal{A} = \sum_{k=1}^{n} Y_k$ and $\mathcal{B} = (2\sum_{k=1}^{n} Y_k^2 + 4p\sum a_k^4)^{1/2}$. It follows from Theorem 4 of [16] that, for any $\lambda \in R$,

$$E\exp\left\{\lambda\mathcal{A} - \frac{\lambda^2}{2}\mathcal{B}^2\right\}$$

$$= \exp\left\{-2\lambda^2 p \sum a_k^4\right\} E \exp\left\{\sum_{k=1}^{n}(\lambda Y_k - \lambda^2 Y_k^2)\right\}$$

$$\leq \exp\left\{-2\lambda^2 p \sum a_k^4\right\}[E\exp\{\lambda Y_1 - \lambda^2 Y_1^2\}]^n$$

$$\leq \exp\left\{-2\lambda^2 p \sum a_k^4\right\}[1 + E(\lambda Y_1 I(\lambda Y_1 \geq -1/2))]^n$$

$$= \exp\left\{-2\lambda^2 p \sum a_k^4\right\}[1 - E(\lambda Y_1 I(\lambda Y_1 \leq -1/2))]^n$$

$$\leq \exp\left\{-2\lambda^2 p \sum a_k^4\right\}[1 + 2\lambda^2 EY_1^2]^n$$

$$\leq \exp\left\{-2\lambda^2 p \sum a_k^4 + 2\lambda^2 n EY_1^2\right\}$$

$$= \exp\left\{-2\lambda^2 p \sum a_k^4 + 2\lambda^2 p \sum(a_k^2 - 1)^2\right\} \leq 1,$$

where we have used the inequality $e^{x-x^2} \leq 1 + xI(x \geq -1/2)$. This yields that two random variables $\mathcal{A}$ and $\mathcal{B} > 0$ satisfy condition (1.4) in [12]. Now, by noting that $(E\mathcal{B})^2 \leq E\mathcal{B}^2 \leq 6p\sum a_k^4$ and applying Theorem 2.1 of [12], we have

$$P\left(V_{1n} \geq 6x\left(\sum_{k=1}^{n}(X_k^2-1)^2 + 5p\sum a_k^4\right)^{1/2}\right)$$

$$\leq P\left(\mathcal{A} \geq \frac{6x}{\sqrt{2}}\sqrt{\mathcal{B}^2 + (E\mathcal{B})^2}\right)$$

(4.4)

$$\leq e^{-6xt/\sqrt{2}} E\exp(t\mathcal{A}/\sqrt{\mathcal{B}^2 + (E\mathcal{B})^2})$$

$$\leq \sqrt{2}e^{-6xt/\sqrt{2}+t^2} \leq \sqrt{2}e^{-4x^2},$$

by letting $t = \sqrt{2}x$. Similarly,

(4.5) $$P\left(-V_{1n} \geq 6x\left(\sum_{k=1}^{n}(X_k^2-1)^2 + 5p\sum a_k^4\right)^{1/2}\right) \leq \sqrt{2}e^{-4x^2}.$$

By virtue of (4.4) and (4.5), we obtain (4.3). The proof of Proposition 2.2 is now complete.

**5. Proof of Proposition 2.3.** Throughout this section, let $\varepsilon_j, 1 \leq j \leq N$, be i.i.d. random variables with $P(\varepsilon_1 = 1) = 1 - P(\varepsilon_1 = 0) = p$, which are also independent of all other random variables, and $B_N = \sum_{j=1}^{N}(\varepsilon_j - p)$. By the inequality $(1+y)^{1/2} \geq 1 + y/2 - y^2$ for any $y \geq -1$, we again have

$$P(S_n \geq x\sqrt{q}V_n)$$

$$= P\left(S_n \geq x\sqrt{nq}\left(1 + \frac{V_n^2 - n}{n}\right)^{1/2}\right)$$

$$\leq P\left(S_n \geq x\sqrt{nq}\left(1 + \frac{V_n^2 - n}{2n} - \frac{(V_n^2-n)^2}{n^2}\right)\right)$$

(5.1)

$$= P\left(\sum \varepsilon_k a_k \geq x\sqrt{nq}\left(1 + \frac{\sum \varepsilon_k(a_k^2-1)}{2n}\right.\right.$$

$$\left.\left. - \frac{(\sum \varepsilon_k(a_k^2-1))^2}{n^2}\right)\Big|B_N = 0\right)$$

$$= P\left(\sum(\varepsilon_k - p)g_k + \frac{x}{n^2}\sum_{1\leq k\neq j\leq N}\nu_k\nu_j \geq x - h\Big|B_N = 0\right)$$

$$= P(T_N + \Lambda_N \geq x - h|B_N = 0),$$

where $h = xpq\sum(a_k^2-1)^2/n^2$,

$$T_N = \sum(\varepsilon_k - p)g_k, \qquad \Lambda_N = \frac{x}{n^2}\sum_{1\leq k\neq j\leq N}\nu_k\nu_j,$$

where, for all $j = 1, \ldots, N$, $\nu_j = (\varepsilon_j - p)(a_j^2 - 1)$ and

$$g_j = \frac{a_j}{\sqrt{nq}} - \frac{x(a_j^2 - 1)}{2n} + \frac{x(1-2p)}{n^2}\left((a_j^2-1)^2 - \frac{1}{N}\sum(a_k^2-1)^2\right),$$

and where in the proof of (5.1) we have used the facts that $\sum a_k = 0, \sum a_k^2 = N$ and

$$(\varepsilon_k - p)^2 = \varepsilon_k(1-2p) + p^2 = (\varepsilon_k - p)(1-2p) + pq.$$

We need the following lemmas before proceeding to the proof of Proposition 2.3.

LEMMA 5.1.  *For any random variable $Z$ with $E|Z| < \infty$,*

$$(5.2) \qquad E(Z|B_N = 0) = \frac{1}{B_n(p)} \int_{-\pi\omega_N}^{\pi\omega_N} EZe^{itB_N/\omega_N}\, dt,$$

*where $B_n(p) = 2\pi\omega_N P(B_N = 0)$ and*

$$(5.3) \qquad 1 \leq \sqrt{2\pi}/B_n(p) \leq 1 + \omega_N^{-2}.$$

PROOF.  Note that $B_N = \sum_{j=1}^N \varepsilon_j - n$ is an integer and, for any integer $k$,

$$\int_{-\pi}^{\pi} e^{ikt}\, dt = \begin{cases} 2\pi, & \text{if } k = 0, \\ 0, & \text{if } k \neq 0. \end{cases}$$

The proof of (5.2) is now obvious. The estimate for $B_n(p)$ follows from $P(B_N = 0) = \binom{N}{n} p^n q^{N-n}$ and Stirling's formula. $\square$

LEMMA 5.2.  (i) *We have*

$$(5.4) \qquad E\left(\sum_{1 \leq k \neq j \leq N} |\nu_k \nu_j|^{3/2} \Big| B_N = 0\right) \leq An^2 \beta_{3N}^2,$$

$$(5.5) \qquad E\left(\sum_{k=1}^N \Big|\nu_k \sum_{j=1, \neq k}^N \nu_j\Big|^{3/2} \Big| B_N = 0\right) \leq An^2 \beta_{3N}^2,$$

$$(5.6) \qquad E\left(\Big|\sum_{1 \leq k \neq j \leq N} \nu_k \nu_j\Big|^{3/2} \Big| B_N = 0\right) \leq An^2 \beta_{3N}^2.$$

(ii) *If $\eta_k, 1 \leq k \leq N$, are i.i.d. random variables with*

$$P(\eta_k = 1) = 1 - P(\eta_k = 0) = m(t), \qquad 0 \leq m(t) \leq 1,$$

*independent of all other random variables, then*

$$(5.7) \quad E\left(\left|\sum_{1\leq k\neq j\leq N}\eta_k\eta_j\nu_k\nu_j\right|^{3/2}\bigg|B_N=0\right)\leq Am^2(t)n^2\beta_{3N}^2,$$

$$(5.8) \quad E\left(\left|\sum_{1\leq k\neq j\leq N}\eta_k(1-\eta_j)\nu_k\nu_j\right|^{3/2}\bigg|B_N=0\right)\leq Am(t)n^2\beta_{3N}^2.$$

PROOF. The proof of Lemma 5.2 is based on an argument similar to that in Theorem 4 of [16], together with the moment inequalities for i.i.d. random variables and $U$-statistics. The details can be found in [18], on which the present paper is based. □

To introduce the following lemmas, we define

$$f(t)=E(e^{it(T_n+\Lambda_n)}|B_N=0), \qquad f_1(t)=E(e^{itT_n}|B_N=0),$$
$$f_2(t)=E(\Lambda_n e^{itT_n}|B_N=0)$$

and, for $k=1,\ldots,N$,

$$g_k(t,\psi)=E\exp\{i(\varepsilon_k-p)(tg_k+\psi/\omega_N)\}.$$

We also use the notation $\Delta=x\beta_{3N}/\omega_N$.

LEMMA 5.3. *If $|t|\leq (1/128)\Delta^{-1}$, then for any $0\leq m(t)\leq 1$ and for $2\leq x\leq (1/128)\omega_N/\max_k |a_k|$,*

$$(5.9) \quad \begin{aligned}|f(t)|&\leq A(1+|tx|)[m^{-1/2}(t)e^{-m(t)t^2/4}+\omega_N e^{-(1/40)m(t)\omega_N^2}]\\ &\quad +A|t|^{3/2}m(t)\Delta^2+A|t|m^{4/3}(t)\Delta^{4/3}.\end{aligned}$$

PROOF. Define $\{\eta_k,k=1,\ldots,N\}$ as in Lemma 5.2(ii). Furthermore, let

$$T_{1N}^*=\sum\eta_k(\varepsilon_k-p)g_k, \qquad T_{2N}^*=\sum(1-\eta_k)(\varepsilon_k-p)g_k,$$
$$\Lambda_{1N}^*=\frac{x}{n^2}\sum_{1\leq k\neq j\leq N}\eta_k\eta_j\nu_k\nu_j, \qquad \Lambda_{2N}^*=\frac{x}{n^2}\sum_{1\leq k\neq j\leq N}\eta_k(1-\eta_j)\nu_k\nu_j,$$
$$\Lambda_{3N}^*=\frac{x}{n^2}\sum_{1\leq k\neq j\leq N}(1-\eta_k)(1-\eta_j)\nu_k\nu_j.$$

Note that

$$(5.10) \quad T_N+\Lambda_N=T_{1N}^*+T_{2N}^*+\Lambda_{1N}^*+2\Lambda_{2N}^*+\Lambda_{3N}^*.$$

It follows from (5.10), $|e^{it} - 1| \leq |t|$ and $|e^{it} - 1 - it| \leq 2|t|^{3/2}$ that

(5.11)
$$\begin{aligned}
|f(t)| &= |E(e^{it(T^*_{1N}+T^*_{2N}+\Lambda^*_{1N}+2\Lambda^*_{2N}+\Lambda^*_{3N})}|B_N=0)| \\
&\leq |E(e^{it(T^*_{1N}+T^*_{2N}+2\Lambda^*_{2N}+\Lambda^*_{3N})}|B_N=0)| + |t|E(|\Lambda^*_{1N}||B_N=0) \\
&\leq |E(e^{it(T^*_{1N}+T^*_{2N}+\Lambda^*_{3N})}|B_N=0)| \\
&\quad + 2|t||E(\Lambda^*_{2N} e^{it(T^*_{1N}+T^*_{2N}+\Lambda^*_{3N})}|B_N=0)| \\
&\quad + 8|t|^{3/2}E(|\Lambda^*_{2N}|^{3/2}|B_N=0) + |t|E(|\Lambda^*_{1N}||B_N=0) \\
&:= \Xi_1(t,x) + \Xi_2(t,x) + \Xi_3(t,x) + \Xi_4(t,x).
\end{aligned}$$

We first estimate $\Xi_3(t,x)$ and $\Xi_4(t,x)$. By Lemma 5.2(ii), we obtain that

$$E(|\Lambda^*_{2N}|^{3/2}|B_N=0) \leq Ax^{3/2}m(t)n^{-1}\beta^2_{3N} \leq Am(t)\Delta^2$$

and, by Hölder's inequality,

$$E(|\Lambda^*_{1N}||B_N=0) \leq [E(|\Lambda^*_{1N}|^{3/2}|B_N=0)]^{2/3} \leq Am^{4/3}(t)\Delta^{4/3}.$$

These facts yield

(5.12) $\quad \Xi_3(t,x) + \Xi_4(t,x) \leq A|t|^{3/2}m(t)\Delta^2 + A|t|m^{4/3}(t)\Delta^{4/3}.$

Next we estimate $\Xi_1(t,x)$. Write $B^*_{1N} = \sum \eta_k(\varepsilon_k - p)$, $B^*_{2N} = \sum(1-\eta_k)(\varepsilon_k - p)$ and

(5.13) $\qquad\qquad B = \{k : \eta_k = 1\}, \qquad B^c = \{k : \eta_k = 0\}.$

Note that, given $\eta_1, \ldots, \eta_N$,

$$T^*_{1N}, B^*_{1N} \in \sigma\{\varepsilon_k, k \in B\}, \qquad T^*_{2N}, \Lambda^*_{3N} B^*_{2N} \in \sigma\{\varepsilon_k, k \in B^c\}$$

and $B_N = B^*_{1N} + B^*_{2N}$. It follows that $T^*_{1N}$ and $B^*_{1N}$ are independent of $T^*_{2N}, \Lambda^*_{3N}$ and $B^*_{2N}$, given $\eta_1, \ldots, \eta_N$, and hence, by Lemma 5.1,

(5.14)
$$\begin{aligned}
\Xi_1(t,x) &= \frac{1}{B_n(p)} \int_{|\psi| \leq \pi\omega_N} |E\exp\{it(T^*_{1N}+T^*_{2N}+\Lambda^*_{3N}) + i\psi B_N/\omega_N\}|\,d\psi \\
&\leq 2\int_{|\psi|\leq\pi\omega_N} E|E_\eta \exp\{itT^*_{1N} + i\psi B^*_{1N}/\omega_N\}|\,d\psi \\
&= 2\int_{|\psi|\leq\pi\omega_N} \prod E|E_\eta\exp\{i\eta_k(\varepsilon_k-p)(tg_k + \psi/\omega_N)\}|\,d\psi,
\end{aligned}$$

where $E_\eta$ denotes the conditional expectation given $\eta_k, k = 1,\ldots,N$.

Let $\varepsilon_k^*$ be an independent copy of $\varepsilon_k$. By Taylor's expansion of $e^{iz}$,

$$\begin{aligned}
E|E_\eta \exp\{i\eta_k(\varepsilon_k - p)(tg_k + \psi/\omega_N)\}|^2 \\
= E(E_\eta \exp\{i\eta_k(\varepsilon_k - \varepsilon_k^*)(tg_k + \psi/\omega_N)\}) \\
= E\exp\{i\eta_k(\varepsilon_k - \varepsilon_k^*)(tg_k + \psi/\omega_N)\} \\
\leq 1 - (1/2)(tg_k + \psi/\omega_N)^2 E\eta_k^2 E(\varepsilon_k - \varepsilon_k^*)^2 \\
+ (1/6)|tg_k + \psi/\omega_N|^3 E\eta_k^3 E|\varepsilon_k - \varepsilon_k^*|^3 \\
\leq 1 - pqm(t)(tg_k + \psi/\omega_N)^2 + (pq/3)m(t)|tg_k + \psi/\omega_N|^3.
\end{aligned}$$

This, together with $\sum g_k = 0$ and the fact that for $2 \leq x \leq (1/128)\omega_N/\max_k |a_k|$,

(5.15) $\quad \left|pq\sum g_k^2 - 1\right| \leq 2x\beta_{3N}/\omega_N \quad \text{and} \quad \sum pq|g_k|^3 \leq 5\beta_{3N}/\omega_N,$

yields that for $2 \leq x \leq (1/128)\omega_N/\max_k |a_k|, |t| < (1/128)\Delta^{-1}$ and $|\psi| < (3/8)\omega_N$,

$$\begin{aligned}
J(t,\psi) &:= \prod E|E_\eta \exp\{i\eta_k(\varepsilon_k - p)(tg_k + \psi/\omega_N)\}| \\
&\leq \left(\prod E|E_\eta \exp\{i\eta_k(\varepsilon_k - p)(tg_k + \psi/\omega_N)\}|^2\right)^{1/2} \\
&\leq \exp\left\{-\tfrac{1}{2}pqm(t)\sum(tg_k + \psi/\omega_N)^2 + \tfrac{1}{6}pqm(t)\sum|tg_k + \psi/\omega_N|^3\right\} \\
&\leq \exp\left\{-(pq/2)m(t)\sum t^2 g_k^2 - m(t)\psi^2/2\right.
\end{aligned}$$

(5.16) $\quad\quad\quad\quad + (2pq/3)m(t)\sum |tg_k|^3 + (2/3)m(t)|\psi|^3/\omega_N\Big\}$

$$\begin{aligned}
&\leq \exp\left\{-(pq/2)m(t)\sum t^2 g_k^2 + (2pq/3)m(t)\sum |tg_k|^3 - m(t)\psi^2/4\right\} \\
&\leq \exp\left\{-\tfrac{1}{2}m(t)t^2(1 - x\beta_{3N}/\omega_N - (5/3)|t|\beta_{3N}/\omega_N) - m(t)\psi^2/4\right\} \\
&\leq \exp\{-m(t)t^2/4 - m(t)\psi^2/4\}.
\end{aligned}$$

To estimate $J(t,\psi)$ for $(3/8)\omega_N \leq |\psi| \leq \pi\omega_N$, we first note that

(5.17)
$$\begin{aligned}
E|E_\eta \exp\{i\eta_k(\varepsilon_k - p)(tg_k + \psi/\omega_N)\}|^2 \\
= E\exp\{i\eta_k(\varepsilon_k - \varepsilon_k^*)(tg_k + \psi/\omega_N)\} \\
= 1 - 2pq + 2pqE\cos[\eta_k(tg_k + \psi/\omega_N)] \\
= 1 - 2pqm(t) + 2pqm(t)\cos(tg_k + \psi/\omega_N).
\end{aligned}$$

Define $D = \{k: |g_k| \leq 2\Delta\}$ and $D^c = \{k: |g_k| > 2\Delta\}$. It is readily seen that for $k \in D$, $|t| < (1/128)\Delta^{-1}$ and $(3/8)\omega_N \leq |\psi| \leq \pi\omega_N$,

$$\tfrac{23}{64} \leq tg_k + \psi/\omega_N \leq \pi + \tfrac{1}{64} \quad \text{or} \quad -\tfrac{1}{64} - \pi \leq tg_k + \psi/\omega_N \leq -\tfrac{23}{64}$$

and hence $\cos(tg_k + \psi/\omega_N) \leq \cos(23/64) < 0.95$. On the other hand, it follows from (5.15) that for $2 \leq x \leq (1/128)\omega_N/\max_k |a_k|$,

$$4(Npq)^{-1}|D^c| \leq \frac{4x^2\beta_{3N}^2}{\omega_N^2}|D^c| \leq \sum_{k \in D^c} g_k^2$$

$$\leq (pq)^{-1}(1 + 2x\beta_{3N}/\omega_N) \leq 2(pq)^{-1},$$

where $|D^c|$ denotes the number of $D^c$. Thus, $|D^c| \leq N/2$ and $|D| = N - |D^c| \geq N/2$. By virtue of (5.17) and all of the above facts, we obtain that for $|t| < (1/128)\Delta^{-1}$, $(3/8)\omega_N \leq |\psi| \leq \pi\omega_N$ and $2 \leq x \leq (1/128)\omega_N/\max_k |a_k|$,

$$J(t, \psi) \leq \left(\prod_{k \in D} E|E_\eta \exp\{i\eta_k(\varepsilon_k - p)(tg_k + \psi/\omega_N)\}|^2\right)^{1/2}$$

(5.18)
$$\leq \prod_{k \in D} \exp\{-pqm(t)[1 - \cos(tg_k + \psi/\omega_N)]\}$$

$$\leq \exp\{-(1/40)m(t)\omega_N^2\}.$$

Combining (5.14), (5.16) and (5.18), it follows that for $|t| < (1/128)\Delta^{-1}$ and $2 \leq x \leq (1/128)\omega_N/\max_k |a_k|$,

(5.19) $\quad \Xi_1(t, x) \leq Am(t)^{-1/2}e^{-m(t)t^2/4} + A\omega_N e^{-(1/40)m(t)\omega_N^2}.$

Finally, we estimate $\Xi_2(t, x)$. Note that $\Lambda_{2N}^* = \frac{x}{n^2}\sum_{j \in B^c}\nu_j \sum_{k \in B}\nu_k$, where $B$ and $B^c$ are defined in (5.13). Similarly to (5.14), we have

$$\Xi_2(t, x) = \frac{2|t|}{B_n(p)}\int_{|\psi| \leq \pi\omega_N} |E(\Lambda_{2N}^* e^{it(T_{1N}^* + T_{2N}^* + \Lambda_{3N}^*) + i\psi B_N/\omega_N})|\,d\psi$$

$$\leq \frac{4|t|x}{n^2}\int_{|\psi| \leq \pi\omega_N} E\Bigg[\sum_{j \in B^c}\sum_{k \in B} E_\eta|\nu_j|$$

(5.20)
$$\times |E_\eta(\nu_k \exp\{itT_{1N}^* + i\psi B_{1N}^*/\omega_N\})|\Bigg]d\psi$$

$$\leq \frac{4|t|x}{n^2}\int_{|\psi| \leq \pi\omega_N} E\Bigg[\sum_{1 \leq j \neq k \leq N}(1 - \eta_j)\eta_k E|\nu_j|E|\nu_k|\Omega_{jk}(t, \psi)\Bigg]d\psi$$

$$\leq \frac{4|t|xm(t)}{n^2}\sum_{1 \leq j \neq k \leq N} E|\nu_j|E|\nu_k|\int_{|\psi| \leq \pi\omega_N} E\Omega_{jk}(t, \psi)\,d\psi,$$

where

$$\Omega_{jk}(t, \psi) = \prod_{l \neq j,k} |E_\eta \exp\{i\eta_l(\varepsilon_l - p)(tg_l + \psi/\omega_N)\}|.$$

As in the proof of (5.19), with minor modifications, we have that for $|t| < (1/128)\Delta^{-1}$, $2 \leq x \leq (1/128)\omega_N / \max_k |a_k|$ and all $1 \leq j \neq k \leq N$,

$$\int_{|\psi| \leq \pi \omega_N} E\Omega_{jk}(t, \psi) \, d\psi \leq Am(t)^{-1/2} e^{-m(t)t^2/4} + A\omega_N e^{-(1/40)m(t)\omega_N^2}.$$

This, together with (5.20) and the fact that

$$\sum_{1 \leq k \neq j \leq N} E|\nu_j| E|\nu_k| \leq \left(2pq \sum (a_k^2 + 1)\right)^2 = 16\omega_N^4,$$

yields that for $2 \leq x \leq (1/128)\omega_N / \max_k |a_k|$ and $|t| < (1/128)\Delta^{-1}$,

(5.21) $$\Xi_2(t, x) \leq A|tx|(e^{-m(t)t^2/4} + \omega_N e^{-(1/40)m(t)\omega_N^2}).$$

Taking estimates (5.12), (5.19) and (5.21) into (5.11), we obtain (5.9). The proof of Lemma 5.3 is now complete. □

LEMMA 5.4. *Suppose that $2 \leq x \leq (1/128)\omega_N / \max_k |a_k|$. Then, for $|t| \leq (1/128)\Delta^{-1/3}$,*

(5.22) $$|f_1(t) - e^{-t^2/2}| \leq A \min\{|t|, 1\}(\Delta(1 + t^6)e^{-t^2/4} + \omega_N^{-6})$$

*and*

(5.23) $$|f_1(t) - g(t, 0)| \leq A \min\{|t|, 1\}(\Delta^{4/3}(1 + t^6)e^{-t^2/4} + \omega_N^{-6}),$$

*where*

$$g(t, 0) = e^{-t^2/2}\left\{1 + \sum(g_k(t, 0) - 1) + t^2/2\right\}.$$

LEMMA 5.5. *Suppose that $2 \leq x \leq (1/128)\omega_N / \max_k |a_k|$. Then, for $|t| \leq (1/128)\Delta^{-1/3}$,*

(5.24) $$|f_2(t)| \leq A(1 + t^2)\Delta^2(e^{-t^2/4} + \omega_N^{-6}),$$

(5.25) $$|f(t) - f_1(t)| \leq A\Delta^2|t|^{3/2} + A|t|(1 + t^2)\Delta^2(e^{-t^2/4} + \omega_N^{-6}).$$

The proofs of Lemmas 5.4 and 5.5 are omitted. The details can be found in [18], on which the present paper is based.

LEMMA 5.6. *Suppose that $2 \leq x \leq (1/128)\omega_N / \max_k |a_k|$. Then there exists an absolute constant $A$ such that, for all $|y| \leq 4x$,*

$$P(T_N + \Lambda_N \geq y | B_N = 0) \leq (1 - \Phi(y)) + Ax\Delta e^{-y^2/2} + A\Delta^{4/3}.$$

PROOF. Note that Lemmas 5.3–5.5 are similar to Lemmas 10.1–10.3 in [19]. The proof of Lemma 5.6 is similar to Lemma 10.5 of [19] with some routine modifications. We omit the details. □

We are now ready to prove Proposition 2.3. Note that $\max |a_k| \leq \omega_N$,

$$h = xpq\sum(a_k^2 - 1)^2/n^2 \leq x\max|a_k|\beta_{3N}/n \leq \Delta$$

and $|x - h| \leq 2x$. It follows from (5.1) and Lemma 5.6 that

$$\begin{aligned}P(S_n \geq x\sqrt{q}V_n) &\leq P(T_N + \Lambda_N \geq x - h|B_N = 0)\\ &\leq (1 - \Phi(x - h)) + Ax\Delta e^{-(x-h)^2/2} + A\Delta^{4/3}\\ &\leq 1 - \Phi(x) + A(1 + x)\Delta e^{-x^2/2 + x\Delta} + A\Delta^{4/3}\\ &\leq (1 - \Phi(x))(1 + Ax^2\Delta e^{x\Delta}) + A\Delta^{4/3}\\ &\leq (1 - \Phi(x))\exp\{Ax^3\beta_{3N}/\omega_N\} + A(x\beta_{3N}/\omega_N)^{4/3},\end{aligned}$$

where we have used the result

$$\Phi(x) - \Phi(x - h) \leq h\Phi'(x - h) \leq he^{-(x-h)^2/2} \leq \Delta e^{-x^2/2 + x\Delta}.$$

This yields Proposition 2.3.

**Acknowledgments.** The authors would like to thank the referees and the Editor for their valuable comments which have led to this much improved version of the paper.

Z. Hu
Department of Statistics and Finance
University of Science
 and Technology of China
Hefei 230026
China
E-mail: huzs@ustc.edu.cn

J. Robinson
Q. Wang
School of Mathematics and Statistics F07
University of Sydney
New South Wales 2006
Australia
E-mail: johnr@maths.usyd.edu.au
       qiying@maths.usyd.edu.au